\documentclass[11pt]{article}
\usepackage{amssymb,amsfonts,amsmath,latexsym, epsfig,mathrsfs}
\newtheorem{theo}{Theorem}
\newtheorem{defn}[theo]{Definition}

\newtheorem{lem} [theo]{Lemma}
\newtheorem{coro}[theo]{Corollary}

\makeatletter \@addtoreset{equation}{section}
\@addtoreset{theo}{}\makeatother

\def\qed{\hfill \rule{4pt}{7pt}}
\def\pf{\noindent {\it Proof.} }

\begin{document}
\title{$(k,m)$-Catalan Numbers and \\ Hook Length Polynomials for Plane Trees}
\author{Rosena R.X. Du\footnote{Department of Mathematics, East China Normal University, Shanghai, P. R. China, rxdu@math.ecnu.edu.cn.}  and Fu Liu\footnote{Department of
Mathematics, Massachusetts Institute of Technology, Cambridge, MA
02139, USA, fuliu@math.mit.edu.}} \maketitle

\noindent {\bf Abstract:} Motivated by a formula of A. Postnikov
relating binary trees, we define the hook length polynomials for
$m$-ary trees and plane forests, and show that these polynomials
have a simple binomial expression. An integer value of this
expression is $C_{k,m}(n)=\frac{1}{mn+1}{(mn+1)k \choose n}$,
which we call the $(k,m)$-Catalan number. For proving the hook
length formulas, we also introduce  a combinatorial family,
$(k,m)$-ary trees, which are counted by the $(k,m)$-Catalan
numbers.



\section{Introduction}
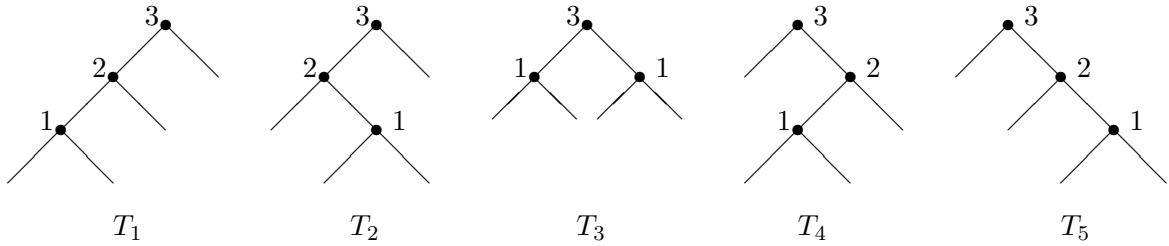
\begin{figure}[h,t]
\begin{center}
\begin{picture}(430,90)
\setlength{\unitlength}{7mm}

\put(2,0){$T_1$} \put(3,4){\line(-1,-1){1}}
\put(3,4){\line(1,-1){1}} \put(3,4){\circle*{0.2}}
\put(2.6,4){$3$}

\put(2,3){\line(-1,-1){1}} \put(2,3){\line(1,-1){1}}
\put(2,3){\circle*{0.2}} \put(1.6,3){$2$}

\put(1,2){\line(-1,-1){1}} \put(1,2){\line(1,-1){1}}
\put(1,2){\circle*{0.2}} \put(0.6,2){$1$}

\put(6.5,0){$T_2$} \put(7,4){\line(-1,-1){1}}
\put(7,4){\line(1,-1){1}} \put(7,4){\circle*{0.2}}
\put(6.6,4){$3$}

\put(6,3){\line(-1,-1){1}} \put(6,3){\line(1,-1){1}}
\put(6,3){\circle*{0.2}} \put(5.6,3){$2$}

\put(7,2){\line(-1,-1){1}} \put(7,2){\line(1,-1){1}}
\put(7,2){\circle*{0.2}} \put(7.3,2){$1$}

\put(10.8,0){$T_3$} \put(11,4){\line(-1,-1){1}}
\put(11,4){\line(1,-1){1}} \put(11,4){\circle*{0.2}}
\put(10.6,4){$3$}

\put(10,3){\line(-1,-1){0.8}} \put(10,3){\line(1,-1){0.8}}
\put(10,3){\circle*{0.2}} \put(9.6,3){$1$}

\put(12,3){\line(-1,-1){0.8}} \put(12,3){\line(1,-1){0.8}}
\put(12,3){\circle*{0.2}} \put(12.3,3){$1$}

\put(15,0){$T_4$} \put(15,4){\line(-1,-1){1}}
\put(15,4){\line(1,-1){1}} \put(15,4){\circle*{0.2}}
\put(15.3,4){$3$}

\put(16,3){\line(-1,-1){1}} \put(16,3){\line(1,-1){1}}
\put(16,3){\circle*{0.2}} \put(16.3,3){$2$}

\put(15,2){\line(-1,-1){1}} \put(15,2){\line(1,-1){1}}
\put(15,2){\circle*{0.2}} \put(14.6,2){$1$}

\put(20,0){$T_5$} \put(19,4){\line(-1,-1){1}}
\put(19,4){\line(1,-1){1}} \put(19,4){\circle*{0.2}}
\put(19.3,4){$3$}

\put(20,3){\line(-1,-1){1}} \put(20,3){\line(1,-1){1}}
\put(20,3){\circle*{0.2}} \put(20.3,3){$2$}

\put(21,2){\line(-1,-1){1}} \put(21,2){\line(1,-1){1}}
\put(21,2){\circle*{0.2}} \put(21.3,2){$1$}
\end{picture}
\caption{Hook lengths for the five binary trees with $3$ internal
vertices.} \label{hook}
\end{center}
\end{figure}

The main result in this paper is originally motivated by seeking a
simple bijective proof of the following identity:
\begin{equation}\label{ihookx=1}
\sum_{T}\frac{n!}{2^n}\prod_{v}\left(1+\frac{1}{h_v}\right)=(n+1)^{n-1},
\end{equation}
where the sum is over all complete binary trees with $n$ internal
vertices, the product is over all internal vertices of $T$, and
$h_v$ is the ``hook length" of $v$ in $T$, namely, the number of
internal vertices in the subtree of $T$ rooted at $v$. As an
example, Figure \ref{hook} shows the five complete binary trees of
three internal vertices, with each internal vertex labelled by its
hook length. In this case, this identity says that
$3+3+4+3+3=(3+1)^2.$

This identity was first derived by Postnikov \cite{apost}, who
also asked for a combinatorial proof of this identity. Chen and
Yang \cite{chenyang} and Seo \cite{Seo} both gave bijective proofs
of it.

Based on (\ref{ihookx=1}), Lascoux conjectured that if we
substitute $x$ for $1$ on the left hand side, we will get the
following polynomial:
\begin{equation}\label{ilas-1}
\sum_{T}\prod_{v}\left(x+\frac{1}{h_v}\right)=\frac{1}{(n+1)!}\prod_{i
= 0}^{n-1} \left((n + 1 + i)x + n +1 -i\right).
\end{equation}

This is equivalent to the following, more suggestive identity:
\begin{equation}\label{ilas-2}
\sum_{T}\prod_{v}\frac{(h_v+1)x + 1 -
h_v}{2h_v}=\frac{1}{n+1}{{(n+1)x} \choose n}.
\end{equation}

We call the left side of (\ref{ilas-2}) the ``hook length
polynomial" of complete binary trees. Note that if we replace $x$
with $k,$ the right hand side of (\ref{ilas-2}) becomes
$\frac{1}{n+1}{{(n+1)k} \choose n}$, which is exactly the number
of complete $k$-ary trees with $n+1$ internal vertices. In fact we
could prove (\ref{ilas-2}) by showing that both the hook length
polynomials for complete binary trees and $k$-ary trees have the
same recurrence relation.

Moreover, we are able to generalize (\ref{ilas-2}) from
enumerating on binary trees to enumerating on $(m+1)$-ary trees,
where $m \geq 0$, and obtain the following identity for hook
length polynomials for $(m+1)$-ary trees:

\begin{equation}\label{ilas-3}
\sum_{T}\prod_{v}\frac{(mh_v+1)x+1-h_{v}}{(m+1)h_v}=\frac{1}{mn+1}{{(mn+1)x}
\choose n}.
\end{equation}

Identity (\ref{ilas-3}) is one of the main results in this paper.
If we replace $x$ with $k,$ the right hand side of (\ref{ilas-3})
becomes
\[\frac{1}{mn+1}{{(mn+1)k} \choose n}.\]
We define this number as the {\it $(k,m)$-Catalan number of order
$n$}. We want to prove (\ref{ilas-3}) using the similar idea as
for (\ref{ilas-2}). Therefore, we need to find some nice
combinatorial interpretation of this number.

In Section 2, we define a combinatorial structure, $(k,m)$-ary
trees, and prove that they are counted by $(k,m)$-Catalan numbers.
In Section 3, we prove (\ref{ilas-3}) by showing that the hook
length polynomials for $(m+1)$-ary trees and $(k,m)$-ary trees
have the same recurrence relation. Some other nice identities
concerning hook lengths for $(m+1)$-ary trees are also studied in
this section. In Section 4, we give similar results for hook
length polynomials for plane forests.

\section{$(k,m)$-ary Trees and $(k,m)$-Catalan numbers}
In this section we define $(k,m)$-ary trees, and prove that they
are counted by $(k,m)$-Catalan numbers.

Let us first review some terminology related to trees. A {\it
tree} is an acyclic connected graph, and a {\it forest} is a graph
such that every connected component is a tree. In this paper we
will assume all the trees are {\it unlabelled plane trees}, i.e.,
rooted trees whose vertices are considered to be
indistinguishable, but the subtrees at any vertex are linearly
ordered. For each vertex of a tree, we say that it is of {\it
degree} $m$ if it has $m$ children, and call vertices of degree
$0$ {\it leaves}. (Note that the definition of degree here is
different from that in graph-theoretic terminology.) Vertices that
are not leaves are called {\it internal vertices}. For any tree
$T$, we use $\mathcal{I}(T)$ to denote the set of internal
vertices of $T$.

A {\it complete binary tree} is a tree each of whose internal
vertices has degree $2$, and a {\it complete $m$-ary tree} is a
tree each of whose internal vertices has degree $m$. Since all of
our trees will be complete, we will frequently just say $m$-ary
tree or binary tree. We will use ${\mathcal T}_{m}(n)$ to denote
the set of all $m$-ary trees with $n$ internal vertices. An
$m$-ary tree with $0$ internal vertices is defined to be a single
vertex.

Let $T$ be a tree with root $r$. For each vertex $v$ of $T$, we
say that $v$ is on {\it level} $j$ if the unique path from $v$ to
$r$ is of length $j$, and the root is said to be on the {\it
level} $0$. We can then define $(k,m)$-ary trees as follows.

\begin{defn}\label{kmtree}
For $k,m \ge 1$ and $n\geq 0,$ a $(k,m)${\rm-ary tree of order}
$n$ is a tree which satisfies the following:
\begin{enumerate}
\item All vertices on even levels have degree $k$. \item All
vertices on odd levels have degree $m$ or $0$, and there are
exactly $n$ vertices of degree $m$.
\end{enumerate}
\end{defn}

\begin{figure}[h,t]
\begin{center}
\begin{picture}(220,100)
\setlength{\unitlength}{7mm} 
\put(5,6){\line(-2,-1){2}} \put(5,6){\line(0,-1){1}}
\put(5,6){\line(2,-1){2}}

\put(3,5){\circle*{0.2}} \put(3,5){\line(-1,-1){1}}
\put(3,5){\line(1,-1){1}} \put(2,4){\line(-1,-2){0.5}}
\put(2,4){\line(0,-1){1}} \put(2,4){\line(1,-2){0.5}}
\put(4,4){\line(-1,-2){0.5}} \put(4,4){\line(0,-1){1}}
\put(4,4){\line(1,-2){0.5}}

\put(7,5){\circle*{0.2}} \put(7,5){\line(-1,-1){1}}
\put(7,5){\line(1,-1){1}} \put(6,4){\line(-1,-2){0.5}}
\put(6,4){\line(0,-1){1}} \put(6,4){\line(1,-2){0.5}}
\put(8,4){\line(-1,-2){0.5}} \put(8,4){\line(0,-1){1}}
\put(8,4){\line(1,-2){0.5}}

\put(6,3){\circle*{0.2}} \put(6,3){\line(-1,-1){1}}
\put(6,3){\line(1,-1){1}} \put(5,2){\line(-1,-2){0.5}}
\put(5,2){\line(0,-1){1}} \put(5,2){\line(1,-2){0.5}}
\put(7,2){\line(-1,-2){0.5}} \put(7,2){\line(0,-1){1}}
\put(7,2){\line(1,-2){0.5}}


%
\end{picture}
\vskip -7mm \caption{A $(3,2)$-ary trees of order $3$.} \label{km}
\end{center}
\end{figure}
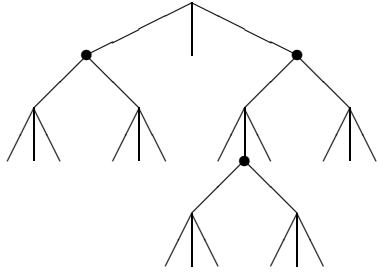

For example, the tree in Figure \ref{km} is a $(3,2)$-ary tree of
order $3$. We will use $\mathcal{T}_{k,m}(n)$ to denote the set of
$(k,m)$-ary trees of order $n$. It is easy to check that there are
a total of $(mn+1)(k+1)$ vertices for each tree $T\in
\mathcal{T}_{k,m}(n)$, with $mn+1$ vertices on even levels (which
are all internal vertices), and $(mn+1)k$ vertices on odd levels
(in which there are $n$ internal vertices and $(mn+1)k-n$ leaves).
We will call the $n$ internal vertices on odd levels {\it crucial
vertices}. Note that a $(k,m)$-ary tree of order $0$ (with $0$
crucial vertices) is a $k$-ary tree with only one internal vertex
(which is the root).

\begin{theo}\label{kmt}
The number of $(k,m)$-ary trees of order $n$ is
\begin{equation}
C_{k,m}(n):=\frac{1}{mn+1}{(mn+1)k \choose n}.
\end{equation}
\end{theo}
\begin{pf} Let $C(z)=\sum_{i=0}^{\infty} c_i z^i$ be the generating function for $(k,m)$-ary
trees, where $c_n$ is the number of $(k,m)$-ary trees of order
$n$.

For $T \in \mathcal{T}_{k,m}(n),$ if we delete the root of $T$, we
will get $k$ subtrees, and each of them can be considered as a
forest of $m$ ordered $(k,m)$-ary trees. Hence we have
\begin{equation}\label{Crecur}
C=(1+zC^m)^k.
\end{equation}
Let $B=zC^m,$ then from (\ref{Crecur}) we have $C = (1+B)^k$ and
$B=z(1+B)^{mk}$. Applying the Lagrange inversion formula (see
\cite[p.42]{EC2}) yields
\begin{equation*}
[z^n]C(z)=\frac{1}{n}[t^{n-1}] \left( \frac{d}{dt} (1+t)^k
\right)(1+t)^{mkn}=\frac{1}{mn+1}{(mn+1)k \choose n}.
\end{equation*}
\end{pf}\qed

\noindent{\it \bf Remarks:}
\begin{enumerate}
\item We can also prove this theorem (as well as Lemma \ref{k2r}
in Section 4) combinatorially, but we omit these proofs for the
sake of brevity. The current proofs of Theorem \ref{kmt} and Lemma
\ref{k2r} were suggested by a referee.

\item It is well known \cite[Ex.~6.19]{EC2} that binary trees with
$n$ internal vertices are counted by the $n$th Catalan number:
\begin{equation}\label{cn}
|\mathcal{T}_{2}(n)|=C(n)=\frac{(2n)!}{(n+1)!n!}=\frac{1}{2n+1}{2n+1
\choose n}=\frac{1}{n}{2n \choose {n-1}}=\frac{1}{n+1}{2n \choose
n},
\end{equation}
and it is also known \cite{Goulden,Klarner,Knuth} that the number
of $m$-ary trees with $n$ internal vertices is
\begin{equation}
|\mathcal{T}_{m}(n)|=C_m(n)=\frac{1}{n} {mn \choose
n-1}=\frac{1}{mn+1} {mn+1 \choose n}.
\end{equation}

For a $(k,m)$-ary tree of order $n$, if $k$ or $m$ equals $1$, we
can simply
 contract the redundant
edges and get the following results:
\begin{eqnarray*}
&\mathcal{T}_{1,2}(n)=\mathcal{T}_{2}(n),\
&\mathcal{T}_{2,1}(n)=\mathcal{T}_{2}(n+1), \\
&\mathcal{T}_{1,m}(n)=\mathcal{T}_{m}(n),\
&\mathcal{T}_{k,1}(n)=\mathcal{T}_{k}(n+1),
\end{eqnarray*}
which coincides with the fact that
\begin{eqnarray}
&C_{1,2}(n)=C(n),\ \
&C_{2,1}(n)=C(n+1), \\
&C_{1,m}(n)=C_m(n),\ \ \ &C_{k,1}(n)=C_k(n+1).
\end{eqnarray}
Hence $(k,m)$-Catalan numbers $C_{k,m}(n)$ can be viewed as a
generalization of $C(n)$ and $C_k(n)$, justifying the choice of
terminology.
\end{enumerate}


\section{Hook length polynomials for $m$-ary trees}

In this section we will study hook length polynomials for $m$-ary
trees. Given an $m$-ary tree $T,$ for any internal vertex $v$,
recall that the {\it hook length} $h_v$ of $v$ is defined as the
number of internal vertices in the subtree rooted at $v$.

We define the {\it hook length polynomial} for an $m$-ary tree $T$
to be
\begin{equation*}
T(x) = \prod_{v \in {\mathcal{I}(T)}}
\frac{((m-1)h_v+1)x+1-h_{v}}{mh_v}.
\end{equation*}
For $n \geq 1$, the {\it $n$th hook length polynomial of $m$-ary
trees} is then defined to be
\begin{equation*}
H_{n,m}(x)=\sum_{T \in \mathcal{T}_m(n)}T(x),
\end{equation*}
and by convention, we set $H_{0,m}(x)=1$. (Recall that
$\mathcal{T}_m(n)$ denotes the set of all $m$-ary trees with $n$
internal vertices.)

The main result of this section is the following identity for the
hook length polynomial of $(m+1)$-ary trees:
\begin{theo}\label{theomhook}
\begin{equation}\label{mhook}
H_{n,m+1}(x) =\frac{1}{mn+1}{(mn+1)x \choose n}.
\end{equation}
\end{theo}

For the example in Figure \ref{hook}, we have
\begin{eqnarray*}
&T_1(x)=T_2(x)=T_4(x)=T_5(x)=\frac{1}{12}x(2x-1)(3x-1),\\
&T_3(x)=\frac{1}{3}x^2(2x-1),
\end{eqnarray*}
and
\begin{equation*}
H_{3,2}(x)=4\cdot\frac{1}{12}x(2x-1)(3x-1)+\frac{1}{3}x^2(2x-1)=\frac{1}{4}{4x
\choose 3},
\end{equation*}
which coincides with equation (\ref{mhook}).

We will prove Theorem \ref{theomhook} by showing that both sides
of (\ref{mhook}) satisfy the same recurrence, as proved in the
following two lemmas.

\begin{lem}\label{bre}
\begin{equation}\label{bref}
H_{n,m+1}(x) = \frac{(mn+1)x+1-n}{(m+1)n}
\sum_{\stackrel{i_1,i_2,\ldots, i_{m+1} \geq 0}{i_1+i_2+ \cdots
+i_{m+1}=n-1}} H_{i_1,m+1}(x) H_{i_2,m+1}(x) \cdots
H_{i_{m+1},m+1}(x).
\end{equation}
\end{lem}
\pf For the sake of convenience, we set
\[h_v(x) =\frac{(mh_v+1)x +1-h_v}{(m+1)h_v}.\]

For any $T \in \mathcal{T}_{m+1}(n),$ if $r$ is the root of $T,$
then $h_r = n$, so $h_r(x) = \frac{(mn+1)x+1-n}{(m+1)n}$. Let
$T_1, T_2, \ldots, T_{m+1}$ be the $(m+1)$-ary trees obtained by
deleting $r$ from $T$, so the total number of internal vertices in
$T_1, T_2, \ldots, T_{m+1}$ is $n-1$, and
$$T(x) =
\frac{(mn+1)x+1-n}{(m+1)n}T_1(x)T_2(x)\cdots T_{m+1}(x).$$ If
$T_j$ has $i_j$ internal vertices for $1 \leq j \leq m+1$, by
summing over all $T \in \mathcal{T}_{m+1}(n)$ we get the desired
result. \qed

\begin{lem}\label{kre}
\begin{equation}\label{kref}
C_{k,m}(n) = \frac{(mn+1)k+1-n}{(m+1)n}
\sum_{\stackrel{i_1,i_2,\ldots, i_{m+1} \geq 0}{i_1+i_2+ \cdots
+i_{m+1}=n-1}} C_{k,m}(i_1)C_{k,m}(i_2)\cdots C_{k,m}(i_{m+1})
\end{equation}
\end{lem}

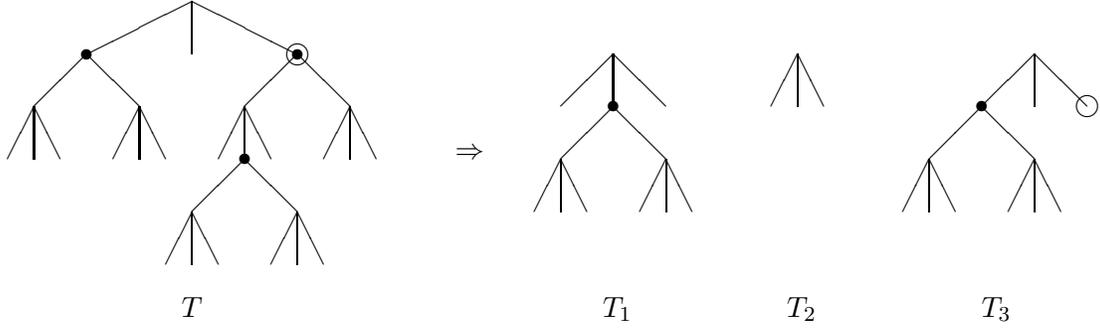
\begin{figure}[h,t]
\begin{center}
\begin{picture}(380,124)
\setlength{\unitlength}{7mm}
\put(2.8,0){$T$}
\put(3,6){\line(-2,-1){2}}
\put(3,6){\line(0,-1){1}}
\put(3,6){\line(2,-1){2}}

\put(1,5){\circle*{0.2}}
\put(1,5){\line(-1,-1){1}}
\put(1,5){\line(1,-1){1}}
\put(0,4){\line(-1,-2){0.5}}
\put(0,4){\line(0,-1){1}}
\put(0,4){\line(1,-2){0.5}}
\put(2,4){\line(-1,-2){0.5}}
\put(2,4){\line(0,-1){1}}
\put(2,4){\line(1,-2){0.5}}

\put(5,5){\circle*{0.2}}
\put(5,5){\circle{0.4}}
\put(5,5){\line(-1,-1){1}}
\put(5,5){\line(1,-1){1}}
\put(4,4){\line(-1,-2){0.5}}
\put(4,4){\line(0,-1){1}}
\put(4,4){\line(1,-2){0.5}}
\put(6,4){\line(-1,-2){0.5}}
\put(6,4){\line(0,-1){1}}
\put(6,4){\line(1,-2){0.5}}

\put(4,3){\circle*{0.2}}
\put(4,3){\line(-1,-1){1}}
\put(4,3){\line(1,-1){1}}
\put(3,2){\line(-1,-2){0.5}}
\put(3,2){\line(0,-1){1}}
\put(3,2){\line(1,-2){0.5}}
\put(5,2){\line(-1,-2){0.5}}
\put(5,2){\line(0,-1){1}}
\put(5,2){\line(1,-2){0.5}}

\put(8,3){$\Rightarrow$}
\put(10.8,0){$T_1$}

\put(11,5){\line(-1,-1){1}}
\put(11,5){\line(0,-1){1}}
\put(11,5){\line(1,-1){1}}
\put(11,4){\circle*{0.2}}
\put(11,4){\line(-1,-1){1}}
\put(11,4){\line(1,-1){1}}
\put(10,3){\line(-1,-2){0.5}}
\put(10,3){\line(0,-1){1}}
\put(10,3){\line(1,-2){0.5}}
\put(12,3){\line(-1,-2){0.5}}
\put(12,3){\line(0,-1){1}}
\put(12,3){\line(1,-2){0.5}}
\put(14.3,0){$T_2$}
\put(14.5,5){\line(-1,-2){0.5}}
\put(14.5,5){\line(0,-1){1}}
\put(14.5,5){\line(1,-2){0.5}}

\put(18,0){$T_3$}
\put(19,5){\line(-1,-1){1}}
\put(19,5){\line(0,-1){1}}
\put(19,5){\line(1,-1){1}}

\put(18,4){\circle*{0.2}}
\put(18,4){\line(-1,-1){1}}
\put(18,4){\line(1,-1){1}}
\put(17,3){\line(-1,-2){0.5}}
\put(17,3){\line(0,-1){1}}
\put(17,3){\line(1,-2){0.5}}
\put(19,3){\line(-1,-2){0.5}}
\put(19,3){\line(0,-1){1}}
\put(19,3){\line(1,-2){0.5}}
\put(20,4){\circle{0.4}}
\end{picture}
\caption{A $(3,2)$-tree of order $3$ with one circled crucial
vertex and the corresponding forest $(T_1,T_2;T_3)$, where $(T_1,
T_2)$ is an ordered pair and $T_3$ has a circled leaf.}
\label{split}
\end{center}
\end{figure}

\pf We will prove this recursion by considering the structure of
$(k,m)$-ary trees. Let $T$ be a tree in $\mathcal{T}_{k,m}(n)$
with one crucial vertex $v$ circled; clearly there are
$nC_{k,m}(n)$ such trees. If we delete the $m$ edges immediately
below $v$, we get a forest of $m$ ordered $(k,m)$-ary trees $T_1,
T_2, \ldots, T_m$ together with a $(k,m)$-ary tree $T_{m+1}$ which
has a leaf circled. Suppose $T_{j}$ has $i_j$ crucial vertices for
all $1 \leq j \leq m+1$, so $i_1+i_2+\cdots+i_{m+1}=n-1$. It is
easy to see that such a split operation is a bijection between
trees in $\mathcal{T}_{k,m}(n)$ with one circled crucial vertex
and the set of $(m+1)$-tuples $(T_1, T_2, \ldots, T_m; T_{m+1})$
such that $T_1, T_2, \ldots, T_m$ are linearly ordered and
$T_{m+1}$ has a circled leaf (see Figure \ref{split}), therefore
we have
\begin{equation}\label{3du1}
\#\{(T_1, T_2, \ldots, T_m;T_{m+1})\}=nC_{k,m}(n).
\end{equation}

Next we will finish the proof by showing that the number of such
tuples equals $n$ times the right hand side of (\ref{kref}). Let
$(T^{\prime}_1, T^{\prime}_2, \ldots, T^{\prime}_m,
T^{\prime}_{m+1})$ be any set of $m+1$ ordered $(k,m)$-ary trees
with a total of $n-1$ crucial vertices such that there is one tree
in the set (not required to be be $T^{\prime}_{m+1}$) which has a
circled leaf. Since a $(k,m)$-ary tree with $i_j$ crucial vertices
has $(mi_j+1)k-i_j$ leaves, there are $(mn+1)k+1-n$ leaves among
all these trees. Hence we have
\begin{equation}\label{3du2}
\#\{(T^{\prime}_1, T^{\prime}_2, \ldots,
T^{\prime}_m, T^{\prime}_{m+1})\}=((mn+1)k+1-n)\cdot
\sum_{\stackrel{i_1,i_2,\ldots, i_{m+1}  \geq 0}{i_1+i_2+ \cdots
+i_{m+1}=n-1}} C_{k,m}(i_1)C_{k,m}(i_2)\cdots C_{k,m}(i_{m+1}).
\end{equation}

On the other hand, for each set $(T^{\prime}_1, T^{\prime}_2,
\ldots, T^{\prime}_m, T^{\prime}_{m+1}),$ we can always set the
tree which has a circled leaf to be the last one to form an
$(m+1)$-tuple $(T_1, T_2, \ldots, T_m; T_{m+1})$, so
\begin{equation}\label{3du3}
\#\{(T^{\prime}_1, T^{\prime}_2, \ldots,
T^{\prime}_m, T^{\prime}_{m+1})\}=(m+1)\#\{(T_1, T_2, \ldots, T_m;
T_{m+1})\}.
\end{equation}

Compare (\ref{3du1}), (\ref{3du2}) and (\ref{3du3}), we get
(\ref{kref}). \qed


\noindent{\it Proof of Theorem \ref{theomhook}:} Since
$H_{n,m+1}(x)$ and $\frac{1}{mn+1}{(mn+1)x \choose n}$ are both
polynomials of $x$ of degree $n$, it is enough to prove that both
sides coincide for positive integer values of $x$. For any
positive integer $k,$ we have
$$C_{k,m}(0) = {k \choose 0} = 1=H_{0,m+1}(k).$$
Moreover, according to Lemma \ref{bre} and Lemma \ref{kre},
$H_{n,m+1}(k)$ and $C_{k,m}(n)$ have the same recurrence relation.
Therefore for any $n \geq 0$, we have
$$H_{n,m+1}(k)=C_{k,m}(n), \forall k \in {\mathbb N}.$$
Hence we have finished the proof. \qed

\begin{coro}\label{lascoux}
\begin{equation}\label{x_1/hv}
\sum_{T \in \mathcal{T}_{m+1}(n)}\prod_{v \in
\mathcal{I}(T)}\left(x+\frac{1}{h_v}\right)=\frac{1}{(mn+1)n!}\prod_{i
= 0}^{n-1} \left((mn + 1 + i)x + mn +1 -mi\right).
\end{equation}
\end{coro}
\pf Define $\varphi_n, E$ and $R_m$ to be the functions that map
polynomials of degree at most $n$ to polynomials of degree at most
$n$, given by $\varphi_n(f(x)) = x^n f(1/x),$ $ E(f(x)) = f(x+1),$
$E^{-1}(f(x)) = f(x-1),$ $R_m(f(x)) = f((m+1)x+m),$ for any
polynomial $f(x)$. Identity (\ref{x_1/hv}) can then be obtained by
applying $E^{-(m-1)} \circ \varphi_n\circ E^{-1} \circ
\varphi_n\circ R_m $ to both sides of (\ref{mhook}).\qed

If we choose the special values $0$ or $m$ for $x$ in
(\ref{x_1/hv}), we will get the following identities.
\begin{coro}
\begin{eqnarray}
n!\sum_{T \in \mathcal{T}_{m+1}(n)} \prod_{v \in {{\mathcal
I}(T)}}
\frac{1}{h_v}&=&\prod_{i=0}^{n-1}(mi+1),\label{mx=0} \\
\sum_{T \in \mathcal{T}_{m+1}(n)} \prod_{v \in
{{\mathcal I}(T)}}
\left(m+\frac{1}{h_v}\right)&=&\frac{(m+1)^{n}(mn+1)^{n-1}}{n!}.\label{mx=1}
\end{eqnarray}
\end{coro}

We want to remark here that equation (\ref{mx=0}) gives the number
of increasing $(m+1)$-ary trees and is also proved by Stanley in
\cite{RSpubthesis}.



\section{Hook length polynomials for plane forests}

We define a {\it plane forest} to be a forest of plane trees which
are linearly ordered, and use $\mathcal{F}(n)$ to denote the set
of plane forests with $n$ vertices. Given a plane forest $F \in
\mathcal{F}(n)$, for any vertex $v$ of $F$, the {\it hook length}
$h_v$ of $v$ is defined as the number of vertices in the subtree
rooted at $v$. (Note that this definition is slightly different
from the definition of hook length for $m$-ary trees. For the
latter one, we count the number of {\it internal} vertices because
we make our $m$-ary trees complete, and the definition in terms of
internal vertices is more natural in this context. But we prefer
to considering complete $m$-ary trees in this paper in order to
make the definition of $(k,m)$-ary trees and many proofs clearer.)

In this section we will study the hook length polynomial for plane
forests.

The {\it hook length polynomial} for a plane forest $F$ is defined
to be
\begin{equation*}
F(x) = \prod_{v \in {V(F)}}\frac{(2h_v-1)x+1-h_{v}}{h_v},
\end{equation*}
where $V(F)$ is the set of vertices of $F.$

The {\it $n$th hook length polynomial of plane forests} is then
defined to be
\begin{equation*}
H_{n}(x)=\sum_{F \in \mathcal{F}(n)}F(x).
\end{equation*}
By convention, we set $H_0(x) = 1.$

The main result of this section is the following identity for the
hook length polynomial of plane forests:
\begin{theo}\label{hookpf}
\begin{equation}\label{hpfe}
H_{n}(x) =\frac{1}{2n+1}{(2n+1)x \choose n}.
\end{equation}
\end{theo}

As in the last section, we will prove (\ref{hpfe}) by showing that
both sides satisfy the same recurrence for integer values of $x$.

\begin{lem}\label{pfr}
\begin{equation}\label{pfre}
 H_{n}(x) = \sum_{i=1}^n \frac{(2i-1)x - (i-1)}{i}
H_{i-1}(x)H_{n-i}(x).
\end{equation}
\end{lem}

\begin{pf}
For $n \ge 1,$ given any plane forest $F \in \mathcal{F}(n),$ we
pick the first plane tree $T$ of $F.$ Suppose $T$ has $i$
vertices, so the rest of the trees in $F$ can be considered as a
plane forest with $n-i$ vertices. Let $V(T)$ be the set of
vertices of $T$ and $r$ the root of $T$, so $h_r = i.$ By removing
$r$ from $T,$ we get another plane forest which has $i-1$ trees.
Therefore,
\begin{eqnarray*}
H_{n}(x) &=& \sum_{i=1}^n \left( \prod_{v \in
{V(T)}}\frac{(2h_v-1)x+1-h_{v}}{h_v} \right) H_{n-i}(x) \\
&=& \sum_{i=1}^n \frac{(2h_r-1)x+1-h_{r}}{h_r}\left(\prod_{v \in
{V(T)}, v \neq r}\frac{(2h_v-1)x+1-h_{v}}{h_v}\right) H_{n-i}(x) \\
&=& \sum_{i=1}^n \frac{(2i-1)x+1-i}{i}\left(\prod_{v \in
{V(T)}, v \neq r}\frac{(2h_v-1)x+1-h_{v}}{h_v}\right) H_{n-i}(x) \\
&=& \sum_{i=1}^n \frac{(2i-1)x+1-i}{i}H_{i-1}(x) H_{n-i}(x).
\end{eqnarray*}
 \qed
\end{pf}

For the right side of (\ref{hpfe}), if we replace $x$ with $k$, we
will get $\frac{1}{2n+1}{(2n+1)k \choose n} = C_{k,2}(n).$ When $k
= 1,$ as we mentioned before $C_{1,2}(n)=C(n)$ is just a Catalan
number, which satisfies the recurrence relation:
$$C(n) = \sum_{i=1}^n C(i-1) C(n-i).$$
The following lemma gives a generalization of this recurrence
relation.

\begin{lem}\label{k2r} For $n \ge 1,$
\begin{equation}\label{k2re}
 C_{k,2}(n) = \sum_{i=1}^n \frac{(2i-1)k - (i-1)}{i}
C_{k,2}(i-1)C_{k,2}(n-i).
\end{equation}
\end{lem}

\begin{pf} As in the proof of Theorem \ref{kmt}, we let $C(z)$ be the generating function for
$(k,2)$-ary trees. And we let $D(z) = z \sum_{i = 0}^{\infty}
\frac{c_i}{i+1}z^i$ be $C(z)$'s antiderivative. By (\ref{Crecur}),
\[C = (1 + zC^2)^k.\]
Differentiating on both sides, we get
\begin{eqnarray*}
C' &=& k(1+zC^2)^{k-1}(C^2 + 2zCC') \\
&=& k \frac{C}{1 + zC^2} (C^2 + 2zCC').
\end{eqnarray*}
One can rearrange this equation to obtain
\[(1-k)C = -(2k-1)(C + zC') + \frac{C'}{C^2}.\]
Then integrating on both sides gives us
\[(1-k)D = -(2k-1)zC - \frac{1}{C} + \alpha, \] where $\alpha$ is a
constant. If we let $z = 0,$ we find that $\alpha = 1.$ The former
equation can be written as
\begin{equation}\label{k2rgf}
\frac{C-1}{z} = (2k-1)C^2 + (1-k)C \frac{D}{z}.
\end{equation}
We complete the proof by observing that (\ref{k2re}) and
(\ref{k2rgf}) are equivalent.
\end{pf}\qed

\noindent{\it Proof of Theorem \ref{hookpf}:} For any positive
integer $k,$ we have:
$$H_0(k) = 1 = \frac{1}{2\cdot 0+1}{(2\cdot 0+1)k \choose 0}=
C_{k,2}(0).$$

According to Lemma \ref{pfr} and Lemma \ref{k2r}, $C_{k,2}(n)$ and
$H_n(k)$ satisfy the same recurrence relation. Therefore, for any
nonnegative integer $n,$
$$\frac{1}{2n+1}{(2n+1)k \choose n}=C_{k,2}(n) = H_n(k), \forall k \in \mathbb{N}.$$
Since both $\frac{1}{2n+1}{(2n+1)x \choose n}$ and $H_n(x)$ are
polynomials of degree $n$ and they agree on infinitely many values
of $x$, we have proved the desired result. \qed

Similar to the proof of Corollary \ref{lascoux}, if we apply $E
\circ \varphi_n \circ E \circ \varphi_n$ to both sides of
(\ref{hpfe}), we get the following identity.
\begin{coro}
\begin{equation}\label{crpf}
\sum_{F \in \mathcal{F}(n)}\prod_{v \in
V(F)}\left(x+\frac{1}{h_v}\right)=\frac{1}{(2n+1)n!}\prod_{i =
0}^{n-1} \left((2n + 1 - i)x + (2n +1-2i)\right).
\end{equation}
\end{coro}

Finally, setting $x = 0$ or $-2,$ we obtain two additional
interesting formulas.
\begin{coro}
\begin{equation}
n!\sum_{F \in \mathcal{F}(n)}\prod_{v \in
V(F)}\frac{1}{h_v}=(2n-1)!!,
\end{equation}
\begin{equation}
\sum_{F \in \mathcal{F}(n)}\prod_{v \in
V(F)}\left(2-\frac{1}{h_v}\right)=\frac{(2n+1)^{n-1}}{n!}.
\end{equation}
\end{coro}


\noindent{\bf Remark:} By comparing (\ref{mhook}) and
(\ref{hpfe}), one may have noticed that the hook length
polynomials for ternary trees and plane forests are the same:
\begin{equation}\label{question}
H_{n,3}(x)=F_n(x).
\end{equation}
In fact, when $m=2$, the right hand side of (\ref{x_1/hv}) and
(\ref{crpf}) differ only in one sign and agree with each other
when $x=0$. This reflects the fact that increasing ternary trees
with $n$ internal vertices and increasing plane forests on $n$
vertices are both counted by $(2n-1)!!$, as proved in
\cite{RSpubthesis} and \cite{chen}, respectively. We wonder if
there is any direct explanation for (\ref{question}).

\vspace{0.5cm}\noindent{\bf Acknowledgments.} The authors would
like to thank Prof.\ Alain Lascoux for his encouragement in
working on this topic, and who conjectured the special case $m=1$
of Corollary \ref{lascoux}. They would also like to thank Prof.\
Richard P. Stanley for helpful discussions and suggestions, and
two referees for their very valuable comments. Most of this work
was done during the first author's visit to the Department of
Mathematics of Massachusetts Institute of Technology, and she
would also like to thank Prof. Stanley for his advice and support.

\end{document}